\documentclass{elsart}

\usepackage{epsfig}

\usepackage{amssymb}
\usepackage{amsmath}
\usepackage{epsfig}
\usepackage{multicol}
\usepackage{color}
\usepackage{cite}
\usepackage{lineno}
\usepackage[dvipdfm,CJKbookmarks,bookmarksopen=true,colorlinks=true,
linkcolor=blue,citecolor=blue,pdfstartview=FitH,pdftitle=title,pdfauthor=lixx]{hyperref}
\newcommand{\beas}{\begin{eqnarray*}}
\newcommand{\eeas}{\end{eqnarray*}}

\begin{document}

\begin{frontmatter}


\title{Adaptive absorbing boundary conditions for Schr\"{o}dinger-type equations:
   application to nonlinear and multi-dimensional problems}
\author{Zhenli Xu$^1$,~Houde Han$^{1,2}$,~Xiaonan Wu$^3$}
\address{1. Department of Mathematics, University of Science and
Technology of China, Hefei, 230026, China.
Email: xuzl@ustc.edu\\
2. Department of Mathematical Sciences, Tsinghua University,
Beijing, 100084, China. E-mail: hhan@math.tsinghua.edu.cn\\
3. Department of Mathematics, Hong Kong Baptist University,
Kowloon, Hong Kong, China. Email: xwu@hkbu.edu.hk}






\begin{abstract}
We propose an adaptive approach in picking the wave-number
parameter of absorbing boundary conditions for
Schr\"{o}dinger-type equations. Based on the Gabor transform which
captures local frequency information in the vicinity of artificial
boundaries, the parameter is determined by an energy-weighted
method and yields a quasi-optimal absorbing boundary conditions.
It is shown that this approach can minimize reflected waves even
when the wave function is composed of waves with different group
velocities. We also extend the split local absorbing boundary
(SLAB) method [Z. Xu and H. Han, {\it Phys. Rev. E}, 74(2006), pp.
037704] to problems in multidimensional nonlinear cases by
coupling the adaptive approach. Numerical examples of nonlinear
Schr\"{o}dinger equations in one- and two dimensions are presented
to demonstrate the properties of the discussed absorbing boundary
conditions.
\end{abstract}

\end{frontmatter}

\section{Introduction}

The numerical solution of partial differential equations on
unbounded domains arises in a large variety of applications in
science and engineering. A typical example we concern in this paper
is the nonlinear Schr\"{o}dinger-type equations in multi-dimensional
space, which describe the gravity waves on deep water in fluid
dynamics, the pulse propagations in optics fibers, and Bose-Einstein
condensations in very low temperature; see Sulem and Sulem
\cite{SS:Book:99} and Agrawal \cite{Agr:BOOK:01} for details.

One principal difficulty to obtain numerical solutions of these
problems is the unboundedness of the physical domain. In order to
overcome this difficulty,  the artificial boundary method
\cite{GIVOLI:BOOK:92,Tsy:ANM:98,Hag:AN:99,Han:FPCAM:06} has been
widely studied in recent decades, with which the original problem
is reduced to an approximate (or equivalent) problem in a
truncated computational domain. The key point of the artificial
boundary method is to construct a ``suitable" artificial boundary
condition on the given artificial boundary for the problem. In
particular, when we consider problems containing wave
propagations, we hope the artificial boundary conditions can
annihilate all the incident waves so as that there is no or minor
reflected waves propagating into the interior domain. These
artificial boundary conditions are also known as absorbing
boundary conditions. For linear problems, many strategies have
been developed to obtain accurate and efficient boundary
conditions, such as \cite{EM:MC:77,Higdon:MC:86} for hyperbolic
wave equations, \cite{HW:JCM:85,YU:JCM:85} for elliptic equations,
and \cite{HR:NM:95,HH:CMA:02-2} for parabolic equations. In the
case of the linear Schr\"{o}dinger equation, there are also
several works
\cite{SY:JCP:97,Arno:VSLID:98,AES:CMS:03,AB:JCP:03,HH:CMS:04,SW:JCP:06,HJW:CMA:05}
developing transparent boundary conditions and studying their
difference approximations and stability. They utilized the
integral transform (Laplace or Fourier transform) or series
expansion method to construct accurate boundary conditions which
are in nonlocal forms. In practical applications the fast
evaluation method \cite{JG:CMA:04} must be developed to discretize
the nonlocal boundary conditions. On the other hand, the authors
in
\cite{Shi:PRB:91,Kus:PRB:92,Di:NFAO:97,FJ:SISC:99,AR:SINA:02,Sze:SINA:04}
constructed absorbing boundary conditions by deriving polynomial
approximations to nonlocal operators in transparent boundary
conditions with Taylor or rational expansions. This class of
boundary conditions  is local, and hence they are easy to
implement.

The treatment of the boundary conditions on the artificial
boundary for nonlinear equations is difficult in general. Hagstrom
and Keller \cite{HK:MC:87} studied some nonlinear elliptic
problems by linearizing the equations. Han et al.
\cite{HWX:JCM:06} and Xu et al. \cite{XHW:CCP:06} respectively
discussed the nonlinear Burgers eqquation and Kardar-Parisi-Zhang
equation. The exact nonlinear artificial conditions have been
obtained by using the Cole-Hopf transformation. For the works
related with the nonlinear Schr\"{o}dinger equations under
consideration, Zheng \cite{Zhe:JCP:06} obtained the transparent
boundary condition using the inverse scattering transform approach
for the cubic nonlinear Schr\"{o}dinger equation in one dimension.
Antonie et al. \cite{ABD:SINA:06} also studied the one-dimensional
cubic nonlinear Schr\"{o}dinger equation and constructed several
nonlinear integro-differential artificial boundary conditions. In
\cite{Sze:CMAME:06,Sze:NM:06}, Szeftel designed absorbing boundary
conditions for one-dimensional nonlinear wave equation by the
potential and the paralinear strategies. Especially, the
one-dimensional nonlinear Schr\"{o}dinger equation was discussed.
The perfectly matched layer (PML) \cite{FL:JOB:05} was also
applied to handling the nonlinear Schr\"{o}dinger equations in
which the nonlinear term can be general.

Recently, Xu and Han \cite{XH:PRE:06} proposed a split local
absorbing boundary (SLAB) method through a time-splitting procedure
to design absorbing boundary conditions for one-dimensional
nonlinear Schr\"{o}dinger equations. The local absorbing boundary
conditions were imposed on the split linear subproblem and yielded a
full scheme by coupling the discretizations for the interior
equation and boundary subproblems. In using local boundary
conditions for the Schr\"{o}dinger-type equations, it is important
to pre-estimate a wave-number parameter (or the group velocity
parameter) of the wave function, as is illustrated in
\cite{XH:PRE:06}, which strongly influences the accuracy of the
boundary condition. In this paper, we present an adaptive parameter
selection approach based on the Gabor transform \cite{Gabor:46} to
capture the wave number near the artificial boundaries in order that
the constructed absorbing boundary conditions can minimize the
reflected wave. In particular, for nonlinear problems, a wave packet
of the nonlinear Schr\"{o}dinger equation will evolve into various
wave packets with different wave numbers. With the Gabor transform,
the boundary conditions can succeed in reflecting the local
structure of the frequency context of the wave. Particular focus of
this paper is to apply the adaptive approach to multi-dimensional
problems with nonlinear terms, in which very few boundary conditions
can work well.

The organization of this paper is the following. In section 2, we
first give a brief overview of absorbing boundary conditions for the
linear Schr\"{o}dinger equation, and then discuss the adaptive
strategy in picking the parameter in boundary conditions.
Two-dimensional boundary conditions for linear problems are also
proposed in this section. In section 3, we are devoted to the
two-dimensional nonlinear Schr\"{o}dinger equation and its numerical
issues in both interior domain and artificial boundaries. Numerical
examples are investigated in section 4.

\section{Absorbing boundary conditions for the linear Schr\"{o}dinger equation}
\subsection{Brief overview of the absorbing boundary conditions}
We shall give a brief overview for local absorbing boundary
conditions of the linear Schr\"{o}dinger equation in one dimension
\begin{equation}i\psi_t=-\psi_{xx}+V\psi,~~~x\in\mathbb{R},~t>0.\label{LSE1}\end{equation}
Set the truncated subdomain $\Omega_i=[x_l,x_r]$ be the
computational domain. Suppose that the potential $V(x)$ is constant
in the exterior domain $\Omega_e=(-\infty,x_l]\cup[x_r,+\infty)$.

Consider the solutions of the Schr\"{o}dinger wave equation
(\ref{LSE1}) in the form of one Fourier mode:
\begin{equation}\psi(x,t)=e^{-i(\omega t-kx)},\end{equation}
where $k$ is the wave number corresponding to space $x$, and
$\omega$ is the time frequency. We have the dual relation between
the space-time $(x,t)$ domain and the wave number-frequency
$(k,\omega)$ domain: $k\leftrightarrow -i\frac{\partial}{\partial
x}$, and $\omega\leftrightarrow i\frac{\partial}{\partial t}$. Using
this duality, we can transform Schr\"{o}dinger wave equation
(\ref{LSE1}) into the Fourier domain resulting in a dispersion
relation to the equation:
\begin{equation}k^2=\omega-V.\label{dis1}\end{equation}

Under the framework of Engquist and Majda approach \cite{EM:MC:77},
solving (\ref{dis1}) in terms of the wave number $k$ gives
\begin{equation}k=\pm\sqrt{\omega-V},\label{pb1}\end{equation}
where the plus sign corresponds to waves moving to the positive $x$
direction, while the minus sign indicates wave motions in opposite
direction. The exact transformation of (\ref{pb1}) to physical space
is nonlocal in time so that one has to save all history data in
memory in order to perform numerical calculations. An effective
substitution is to approximate the square root through a rational
polynomial.

Let us consider the right exterior domain and obtain boundary
condition at $x=x_r$; that is, the plus sign is taken in
(\ref{pb1}). Similar procedure can be performed in the left exterior
domain. As in \cite{AR:SINA:02}, we denote the absorbing boundary
condition by ABC$(j_1,~j_2)$ for that using $(j_1,~j_2)$-Pad\'{e}
approximation, where $j_1,~j_2$ are the degrees of the polynomials
in the numerator and denominator, respectively.

The first absorbing boundary condition is the one developed in
Shibata \cite{Shi:PRB:91}. The author used a linear interpolation
to approximate the square root in (\ref{pb1}) through imposing two
adjustable parameters which were positive and called the kinetic
energy parameters related to the group velocities of the wave
function \cite{FJ:SISC:99}, that is,
\begin{equation}\sqrt{\omega-V}=\frac{1}{\alpha_1+\alpha_2}(\omega-V)+\frac{\alpha_1\alpha_2}{\alpha_1+\alpha_2}.\end{equation}
Then using the dual relations to transform back into physical space
yields an absorbing boundary condition
\begin{equation}i(\alpha_1+\alpha_2)\psi_x+(i\psi_t-V\psi+\alpha_1\alpha_2\psi)=0.\end{equation}

Kuska \cite{Kus:PRB:92} used a $(1,1)$-Pad\'{e} approximation to
$k^2$ centered at a positive constant $k=k_0$
\begin{equation}k^2=k_0^2\frac{-3k+k_0}{k-3k_0}+O((k-k_0)^3),\end{equation}
and then obtained a second absorbing boundary condition ABC(1,1),
\begin{equation}-\psi_{xt}+i(3k_0^2-V)\psi_x+(k_0^3-3k_0V)\psi+ i3k_0\psi_t=0\label{abc1}\end{equation}
after transforming back into physical space through the dual
transform. Here the range of validity of the Pad\'{e} polynomial is
the positive part $k>0$. Alonso-Mallo and Reguera \cite{AR:SINA:02}
developed a class of absorbing boundary conditions including
ABC(2,1), ABC(3,2) and ABC(2,0) and absorbing boundary conditions
discussed above. Fevens and Jiang \cite{FJ:SISC:99} developed a
distinct method to construct absorbing boundary conditions. The
authors used the group velocity $C=\frac{\partial\omega}{\partial
k}=2k$ to design a differential equation as absorbing boundary
condition, which can absorb waves with certain group velocities
$C_l,~l=1,\cdots,p$,
\begin{equation}\prod_{l=1}^p\left(i\partial_x+\frac{C_l}{2}\right)\psi=0.\label{FJ}\end{equation}

If we substitute the temporal derivative into ABC(1,0) with the
original equation $i\psi_t=-\psi_{xx}+V\psi$, then we obtain
\begin{equation}
(i\partial_x+\alpha_1)(i\partial_x+\alpha_2)\psi=0.
\end{equation}
It is a special case of Fevens and Jiang's formula (\ref{FJ}) with
$C_1=2\alpha_1$ and $C_2=2\alpha_2$. We use the original equation
again to replace the temporal derivative terms in (\ref{abc1}) and
get
\begin{equation}\left(i\frac{\partial}{\partial x}+ k_0\right)^3\psi=0,\end{equation}
which is also a special case of Eq. (\ref{FJ}) for $p=3$ and group
velocities $C_1=C_2=C_3=2k_0$.

\subsection{Weighted wave-number parameter based on Gabor
transform} In the above absorbing boundary conditions, the authors all
imposed parameters in the formulae with different meanings.
Therefore, perhaps one of the most important issues  is how to
pick suitable parameters such that they can minimize the
reflection of the wave.  Noticing that the relation between the
group velocity $C$ and wave number $k$ is
\begin{equation}C=\frac{\partial\omega}{\partial k}=2k,\end{equation}
we need only calculate one of them.

For the initial wave composed of waves with different group
velocities, they shall evolve into different wave packets. Each of
them has an unchanged group velocity. These wave packets hit the
artificial boundary separately. Therefore, in a general physical
insight, if only we pre-estimate one component of group velocities
which is a function of time, the boundary condition can well
annihilate the reflected wave. Let us consider the third-order
boundary condition ABC(1,1) given in section 2.1 as example to
introduce our idea, in which only one parameter $k_0$ need to be
pre-estimated. Similarly, for convenience, the discussion is
focused on the right boundary.

It is important that we must estimate the parameter in the
frequency domain. Note that the wave function at time $t$ can be
expressed in terms of a Fourier series and a single Fourier mode
is essentially a plane wave. A general strategy suggested in
Fevens and Jiang \cite{FJ:SISC:99} to pick the wave-number
parameter $k_0$, which is a function of time $t$, is to use a
Fourier series expansion of the physical variable in space, and
then take one of the positive components so that its Fourier mode
is dominant. The Fourier transform presents the frequency
information of the wave over the whole interior domain. However,
in our situations to construct absorbing boundary conditions, we
are interested in the frequency content of the wave in the
vicinity of the artificial boundary. So it is necessary to obtain
the local structure of the wave in the frequency domain. One
approach is to replace the Fourier transform with the Gabor
transform which is also known as a windowed Fourier transform. In
the frequency domain with the Gabor transform, we have
\begin{equation}\hat{\psi}(k,t)=\int_{x_l}^{x_r}W(x)\psi(x,t)e^{-ikx}
=\int_{x_r-b}^{x_r}\psi(x,t)e^{-ikx},\end{equation} where the
window function is
\begin{equation}
W(x)=\left\{\begin{array}{ll}1,~~~~x\in[x_r-b, x_r],\\0,
~~~~\mathrm{otherwise},\end{array}\right.
\end{equation}
and $b$ is the window width. Then one choice for $k_0$ is take the
frequency such that its spectrum is the maximum; that is,
\begin{equation}|\hat{\psi}(k_0,t)|=\sup_{k\ge 0}\{|\hat{\psi}(k,t)|\}.
\label{max}\end{equation}

We remark that we can also utilize the time windowed Fourier
transform to approximate temporal frequency information $\omega$,
and then obtain an estimation of the wave number $k_0$ with the
dispersion relation (\ref{pb1}). However, it is clear that the
Gabor transform in time depends on the history data on the
artificial boundaries. Therefore, although the formulae of
absorbing boundary conditions are in local forms, they are
nonlocal in practice.

The formula (\ref{max}) is also not the best choice in many
practical computations. On one hand, this procedure involves many
logical ``if" structures in order to compare the magnitudes of the
Fourier modes, which are not very efficient in calculations in
some computational environments. On the other hand, when two
Fourier modes are both dominant, it is obvious to choose $k_0$ a
medial value of two different wave numbers instead of taking one
of them, in order to minimize the reflection. Therefore, an
improvement is to use a weighted strategy, we call it the
energy-weighted wave-number parameter selection approach, as
follows,
\begin{equation}k_0=\int_0^{\infty}(|\hat{\psi}(k,t)|^pk)dk\left/\int_0^{\infty}|\hat{\psi}(k,t)|^p,\right.\label{weight}\end{equation}
with $p$ a positive real number.

We give the following remarks:

{\bf Remark 1.} The window width $b$ is correlative with the Gibbs
phenomena induced by the discontinuities of the window function.
The narrower $b$ is, the more Gibbs effect. However, if the window
width $b$ is very large, then the obtained parameter cannot
correctly response the frequency information in the vicinity of
the boundary.

{\bf Remark 2.} When $p=+\infty$, the Eq. (\ref{weight}) is
equivalent to (\ref{max}). However, numerical experiments
illustrate that the absorbing boundary conditions work best when
$p$ is in a suitable intermediate interval. Table 1 suggests $p=4$
is a good choice.

\subsection{Multi-dimensions} Let us consider the extension of previous ABCs which
are local for the linear Schr\"{o}dinger equation in two dimensions:
\begin{equation}
i\psi_t=-(\psi_{xx}+\psi_{yy})+V\psi,
~~~(x,~y)\in\mathbb{R}^2,\label{LSE}
\end{equation}
with the potential $V$ constant. Denote the dual variables to
$(x,y,t)$ by $(\xi,\eta, \omega)$ with the correspondence
$\xi\leftrightarrow -i\frac{\partial}{\partial
x},~\eta\leftrightarrow -i\frac{\partial}{\partial y},$ and
$\omega\leftrightarrow i\frac{\partial}{\partial t}$. Then the
related dispersion relation to Eq. (\ref{LSE}) gives
\begin{equation}\xi^2+\eta^2=\omega-V.\label{dis}\end{equation}

We truncate the unbounded domain to get a computational domain
$[0~,L]^2$. Without loss of generality, consider the east boundary
$\Gamma_e=\{(x,y)|x=L,~0\leq y\leq L\}$ which corresponds to the
positive branch to $\xi$ of the dispersion relation (\ref{dis}) as
follows,
\begin{equation}\xi=\sqrt{\omega-V-\eta^2}. \label{pb}\end{equation}
With the same procedure as that used in one dimensional case, we can
get the similar ABCs as in section 2.1.  We consider the
$(1,1)-$Pad\'{e} approximation to the square $\xi^2$ in the
dispersion relation centered as a positive constant $\xi=\xi_0$, and
obtain an approximation to (\ref{pb})
\begin{equation}(\eta^2-\omega+V-3\xi_0^2)\xi+\xi_0^3-3\xi_0(\eta^2-\omega+V)=0, \label{pbap}\end{equation}
which is first order in $\xi$. Here the range of validity of the
Pad\'{e} polynomial is the positive part $\xi>0$; see Kuska
\cite{Kus:PRB:92}. Transforming (\ref{pbap}) back into the physical
space through the dual relations, we have an ABC on the right
boundary of the form:
\begin{equation}\Gamma_e:~~~~i\psi_{xyy}-\psi_{xt}+i(3\xi_0^2-V)\psi_x+(\xi_0^3-3\xi_0V)\psi+3\xi_0\psi_{yy}+i\xi_0\psi_t=0. \label{east}\end{equation}
Absorbing boundary conditions on the west, north and south
boundaries can also be obtained through using $(1,~1)$-Pad\'{e}
approximations to $\xi^2$ centered at $-\xi_0$, to $\eta^2$ centered
at $\eta_0$, and to $\eta^2$ centered at $-\eta_0$, respectively,
which are
\begin{equation}\Gamma_w:~~~~i\psi_{xyy}-\psi_{xt}+i(3\xi_0^2-V)\psi_x-(\xi_0^3-3\xi_0V)\psi-3\xi_0\psi_{yy}-i\xi_0\psi_t=0,\end{equation}
\begin{equation}\Gamma_n:~~~~i\psi_{xxy}-\psi_{yt}+i(3\eta_0^2-V)\psi_y+(\eta_0^3-3\eta_0V)\psi+3\eta_0\psi_{xx}+i\psi_t=0,\end{equation}
\begin{equation}\Gamma_s:~~~~i\psi_{xxy}-\psi_{yt}+i(3\eta_0^2-V)\psi_y-(\eta_0^3-3\eta_0V)\psi-3\eta_0\psi_{xx}-i\psi_t=0,\end{equation}
with $\eta_0$ a positive constant as $\xi_0$.

Now let us look at the formula at the north east corner
$(x,y)=(L,L)$. We can also approximate the two dimensional
dispersion relation (\ref{dis}) in the quarter
$\{(\xi,\eta):~\xi>0,~\eta>0\}$ using $(1,1)-$Pad\'{e} to both
$\xi^2$ and $\eta^2$ with the corresponding centered point
$(\xi_0,\eta_0)$ to obtain
\begin{equation}\xi_0^2\frac{-3\xi+\xi_0}{\xi-3\xi_0}+\eta_0^2\frac{-3\eta+\eta_0}{\eta-3\eta_0}=\omega-V\end{equation}
Then after multiplying $(\xi-3\xi_0)(\eta-3\eta_0)$ in both sides,
we have,
\begin{eqnarray}&&-\omega\xi\eta+3\xi_0\omega\eta+3\eta_0\omega\xi+(V-3\xi_0^2-3\eta_0^2)\xi\eta-9\xi_0\eta_0\omega
+(\xi_0^3+9\xi_0\eta_0^2-3\xi_0V)\eta\nonumber\\
&&~~~~~~~~+(\eta_0^3+9\xi_0^2\eta_0-3\xi_0V)\xi
+(9\xi_0\eta_0V-3\xi_0^3\eta_0-3\eta_0^3\xi_0)=0.\end{eqnarray} Then
performing the inverse transform to physical space yields the
ABC(1,1) at the corner point,
\begin{eqnarray}&&i\psi_{xyt}+3\xi_0\psi_{yt}+3\eta_0\psi_{xt}+(3\xi_0^2+3\eta_0^2-V)\psi_{xy}-9i\xi_0\eta_0\psi_t
-i(\xi_0^3+9\xi_0\eta_0^2-3\xi_0V)\psi_y\nonumber\\
&&~~~~~~~~-i(\eta_0^3+9\xi_0^2\eta_0-3\xi_0V)\psi_x
+(9\xi_0\eta_0V-3\xi_0^3\eta_0-3\eta_0^3\xi_0)\psi=0.\label{corner}
\end{eqnarray}

Extension of the adaptive parameter selection for one-dimensional
version to multidimensional cases is straightforward. We note that
a multidimensional problem can be split into  a series of
one-dimensional ones. Thus we can obtain the estimation of the
parameters at every boundary grid points by a
dimension-by-dimension procedure. For example, in order to compute
the wave number on the east boundary, we have the Gabor transform
in $x$ direction:
\begin{equation}\tilde{\psi}(\xi,y,t)=\int_{L-b(y)}^{L}\psi(x,y,t)e^{-i\xi x}dx,\end{equation}
where the window length is a function of $y$. The parameter
$\xi_0(y)$ can then be determined by using the method in section
2.2.

\section{Nonlinear Schr\"{o}dinger equations}

We now consider the nonlinear Schr\"{o}dinger equation in two
dimensions as follows,
\begin{equation}
i\psi_t(x,y,t)=-(\psi_{xx}+\psi_{yy})+f(|\psi|^2)\psi+V(x,y,t)\psi,~~~(x,y)\in\mathbb{R}^2.\label{twod}
\end{equation}
We shall extend the previous work of the split local absorbing
boundary (SLAB) method \cite{XH:PRE:06} in one-dimensional version
to solving the two-dimensional case. Denote the approximation of
$\psi$ on the grid point $(x_\iota ,y_j,t^n)$ by $\psi_{\iota j}^n$
for $0\leq\iota\leq I$ and $0\leq j\leq J$, with $x_\iota
=\iota\Delta x$, $y_j=j\Delta y$, $t^n=n\Delta t$, and $x_I=y_J=L$.
Let us first describe the finite difference scheme for the Eq.
(\ref{twod}) in the interior domain $(0,L)^2$, which will be
connected with the discretization on the artificial boundaries.

\subsection{Semi-implicit interior scheme}

In our previous work in one dimension \cite{XH:PRE:06}, the
full-implicit Crank-Nicholson scheme, which is unconditionally
stable, was used. However, one has to solve the nonlinear algebraic
system iteratively at each time step. It is time consuming in
particular for the two dimensional case. In order to avoid the
iterative process, we use the following semi-implicit scheme
\cite{SSP:CPAM:84}, which was shown efficient and robust in
comparison with various difference schemes for solving nonlinear
Schr\"{o}dinger equations \cite{CJS:JCP:99},
\begin{equation}
i\frac{\psi_{\iota j}^{n+1}-\psi_{\iota j}^n}{\Delta
t}=-(D_x^+D_x^-+D_y^+D_y^-)\frac{\psi_{\iota j}^{n+1}+\psi_{\iota
j}^n}{2} +[\frac{3}{2}f(|u_{\iota j}^n|^2)-\frac{1}{2}f(|u_{\iota
j}^{n-1}|^2)+V_{\iota j}]\frac{\psi_{\iota j}^{n+1}+\psi_{\iota
j}^n}{2}, \label{scheme}
\end{equation}
where $D^+$ and $D^-$ represent the forward and backward
differences, respectively. This is a five-points scheme and its
truncation error is order $O(\Delta t^2+\Delta x^2+\Delta y^2)$ as
that of the Crank-Nicholson scheme. However, since the nonlinear term
is approximated by the known variables through an extrapolation
formula, we  need only to solve a linear algebraic system at each
time step.

\subsection{Numerical approximation on the artificial boundary}

We have obtained the discrete scheme by formula (\ref{scheme}) in
the interior point $(x_\iota ,y_j)$ for $1\leq \iota\leq I-1$ and
$1\leq j\leq J-1$. Now we concentrate on the boundary conditions, in
which we shall perform the local time-splitting procedure. The basic
idea of the SLAB method is to split the original equation in several
subproblems which are easy to be handled, and then solve them
alternatively in a small time step $\Delta t$. Consider a standard
splitting for Eq. (\ref{twod}) in the vicinity of the artificial
boundary to a nonlinear subproblem
\begin{equation}
i\psi_t=f(|\psi|^2)\psi,\label{ns}
\end{equation}
and a linear subproblem
\begin{equation}
i\psi_t=-(\psi_{xx}+\psi_{yy})+V\psi.\label{ls}
\end{equation}
We carry out the splitting on boundary points $\{x_\alpha,y_\beta\}$
for
$$\alpha\in\{0,1,I-1,I\},~~~\mathrm{and}~~~\beta\in\{0,1,J-1,J\}.$$
Following \cite{BJM:SISC:03}, in the nonlinear step, we have an
approximate solver for explicitly discretizing the ODE (\ref{ns})
\begin{equation}
\psi_{\alpha,\beta}^*=e^{-if(|\psi_{\alpha,\beta}^n|^2)\Delta
t}\psi_{\alpha,\beta}^n,
\end{equation}
which keeps $\psi$ invariant, and does not require any boundary
condition. Noting the next step for the time-splitting procedure is
to integrate a linear subproblem (\ref{ls}), we impose here the
local absorbing boundary condition discussed in Section 2. For
example, using formulae (\ref{east}), (\ref{corner}) and their
corresponding formulae on every boundaries and corners, we obtain
the full scheme of the problem by approximating them with finite
difference expressions. Here the discrete forms of the terms in the
east boundary condition (\ref{east}) are
\begin{eqnarray}
&\psi_x=D_x^-\frac{\psi_{I,j}^{n+1}+\psi_{I,j}^*}{2},&\psi=S_x^-\frac{\psi_{I,j}^{n+1}+\psi_{I,j}^*}{2},\\
&\psi_{xt}=D_x^-\frac{\psi_{I,j}^{n+1}-\psi_{I,j}^*}{\Delta t},&\psi_t=S_x^-\frac{\psi_{I,j}^{n+1}-\psi_{I,j}^*}{\Delta t},\\
&\psi_{xyy}=D_x^-D_y^+D_y^-\frac{\psi_{I,j}^{n+1}+\psi_{I,j}^*}{2},
&\psi_{yy}=S_x^-D_y^+D_y^-\frac{\psi_{I,j}^{n+1}+\psi_{I,j}^*}{2},
\end{eqnarray}
with $S^-$ the backward sum, for example,
$$S_x^-\psi_{I,j}^*=\frac{1}{2}(\psi_{I-1,j}^*+\psi_{I,j}^*);$$
and the discrete forms of the terms in the corner boundary condition
(\ref{corner}) are
\begin{eqnarray}
&\psi_{xyt}=D_x^-D_y^-\frac{\psi_{I,J}^{n+1}-\psi_{I,J}^*}{\Delta
t},
&\psi_{xy}=D_x^-D_y^-\frac{\psi_{I,J}^{n+1}+\psi_{I,J}^*}{2},\\
&\psi_{yt}=S_x^-D_y^-\frac{\psi_{I,J}^{n+1}-\psi_{I,J}^*}{\Delta t},
&\psi_y=S_x^-D_y^-\frac{\psi_{I,J}^{n+1}+\psi_{I,J}^*}{2},\\
&\psi_{xt}=D_x^-S_y^-\frac{\psi_{I,J}^{n+1}-\psi_{I,J}^*}{\Delta t},
&\psi_x=D_x^-S_y^-\frac{\psi_{I,J}^{n+1}+\psi_{I,J}^*}{2},\\
&\psi_t=S_x^-S_y^-\frac{\psi_{I,J}^{n+1}-\psi_{I,J}^*}{\Delta
t},&\psi=S_x^-S_y^-\frac{\psi_{I,J}^{n+1}+\psi_{I,J}^*}{2}.
\end{eqnarray}
Similar discretizations can be used for the other three boundaries
and the other three corners. Thus we obtain the full-discrete scheme
for the nonlinear Schr\"{o}dinger equation (\ref{twod}) in two
dimensions, which yields a linear algebraic system at each time
steps.

{\bf Remark 3.} Near the artificial boundaries, the truncation error
of accuracy is $(\Delta x^2+\Delta y^2+\Delta t)$ because we only
adopt the first-order splitting. In order to improve the accuracy of
time-splitting to higher order, such as using the Strang splitting
\cite{Str:SINA:68}, we will obtain a nonlinear algebraic system
which have to be solved through the iterate approach as discussed in
Ref. \cite{XH:PRE:06}.

\section{Numerical examples}

We test our absorbing boundary conditions given in the previous
sections for the nonlinear Schr\"{o}dinger equation. In
particular, we test the strategy of adaptive parameter selection
in the one-dimensional case. Based on its outstanding performance
in one dimension, two-dimensional example of extensions is also
given.

{\bf Example 1.} We are going to test the performance of the
adaptive parameter selection for absorbing boundary conditions by
solving the cubic nonlinear Schr\"{o}dinger equation in one
dimension
\begin{equation}
i\psi_t=-\psi_{xx}+g|\psi|^2\psi+V\psi, \label{oned}
\end{equation}
where $g$ is a real constant and $V\equiv 0$. If $g$ is positive,
the equation represents repulsive interactions. If $g$ is negative,
the equation represents attractive interactions, and admits bright
soliton solution
\begin{equation}
\psi(x,t)=A\sqrt{\frac{-2}{g}}\sec(Ax-2ABt)e^{iBx+6(A^2-B^2)t}
\end{equation}
with $A,B$ real parameters related to the amplitude and velocity of
the soliton. The numerical scheme is the 1D reduction of 2D version
described in Section 3. The first example we consider the case of
$g=-2$ and the initial condition
\begin{equation}
\psi(x,t)=\sec(x-10)e^{2i(x-10)}+\sec(x-30)e^{5i(x-30)}.
\end{equation}
It represents two solitons with amplitude 1, located at two isolated
centers $x=10$ and 30 respectively, propagate to the right. Their
propagating velocities are double of their wave numbers; that is, 4
and 10, respectively. We compute the solution up to $t^n=10$ in
interval $[0,L]$ for $L=40$. As in
\cite{Kus:PRB:92,FJ:SISC:99,XH:PRE:06}, to see the influence of
parameter $k_0$, we evaluate the effectiveness of absorbing boundary
conditions by calculating the reflection ratios as follows,
\begin{equation}r=\sum_{j=0}^I|\psi_j^n|^2/\sum_{j=0}^I|\psi_j^0|^2.\end{equation}
The ratio $r$ is handy in the measurement of the quality of the
ABC. For example, $r=0$ reflects that the solitons have passed
through the boundary completely; whereas $r=1$ indicates the waves
are completely reflected into the interior domain by the
artificial boundaries. At the left boundary $x=0$, we set $k_0=0$
since there is no left-going wave. We also hope the reflected
waves from the right boundary can also be reflected by the left
boundary, therefore the reflection ratios of absorbing boundary
conditions at the right boundary can be correctly calculated. We
show numerical results in Table 1 for different $p$ in
(\ref{weight}) and different transforms, in which we also
illustrate $L_1$-errors defined by
\begin{equation}E_1=\frac{1}{I+1}\sum_{j=0}^I|\psi_j^n-\psi(x_j,t^n)|.\end{equation} Here and
hereafter, the time steps $\Delta t$ are taken to be $\Delta
t=\Delta x^2$. It is not the restriction of stability, but the
requirement for compensating the accuracy since we just use the
first-order splitting on the artificial boundaries. For the Gabor
transform to pre-estimate the parameter, the window lengths are
set to $b=L/4$. We also compare the results in Table 2 without the
adaptive parameter selection but fixing the parameter $k_0=2,~3.5$
and 5, respectively. It is observed that the weighted wave-number
parameter method can well improve numerical accuracy.

There are two time phases in the process for the two solitons hit
the right boundary separately. The first phase is for $t\in
[0.5,~1.5]$ when the first wave with the wave number 5 transmits
the boundary; while the second phase is for $t\in[6,~8]$ when the
second wave with the wave number 2 passes through the boundary. In
order to see the resultant wave-number parameter of the methods in
discussion, we illustrate the wave numbers as a function of time
$t$ for two phases in Fig. 1, where we denote the resultant wave
number with Fourier transform and $p$ norm by F$p$ wave number,
the results with Gabor transform and $p$ norm by G$p$ wave number.
We see that the parameters with Gabor transform, especially when
$p=4$, response a better information for the solution, which well
agrees with the results in Table 1.

  In order to see the influence of the window length $b$, we compute
the solution for different lengths in Table 3 with a fixed $p$.
The window length is in direct proportion to the wave number
$b=\beta k_0$, in which $k_0$ takes the value at time $t=t^{n-1}$
for the calculations at $t^n$. We see that it is necessary to
choose a $\beta$ larger than 1.

\begin{center}\begin{figure}
 \epsfig{file=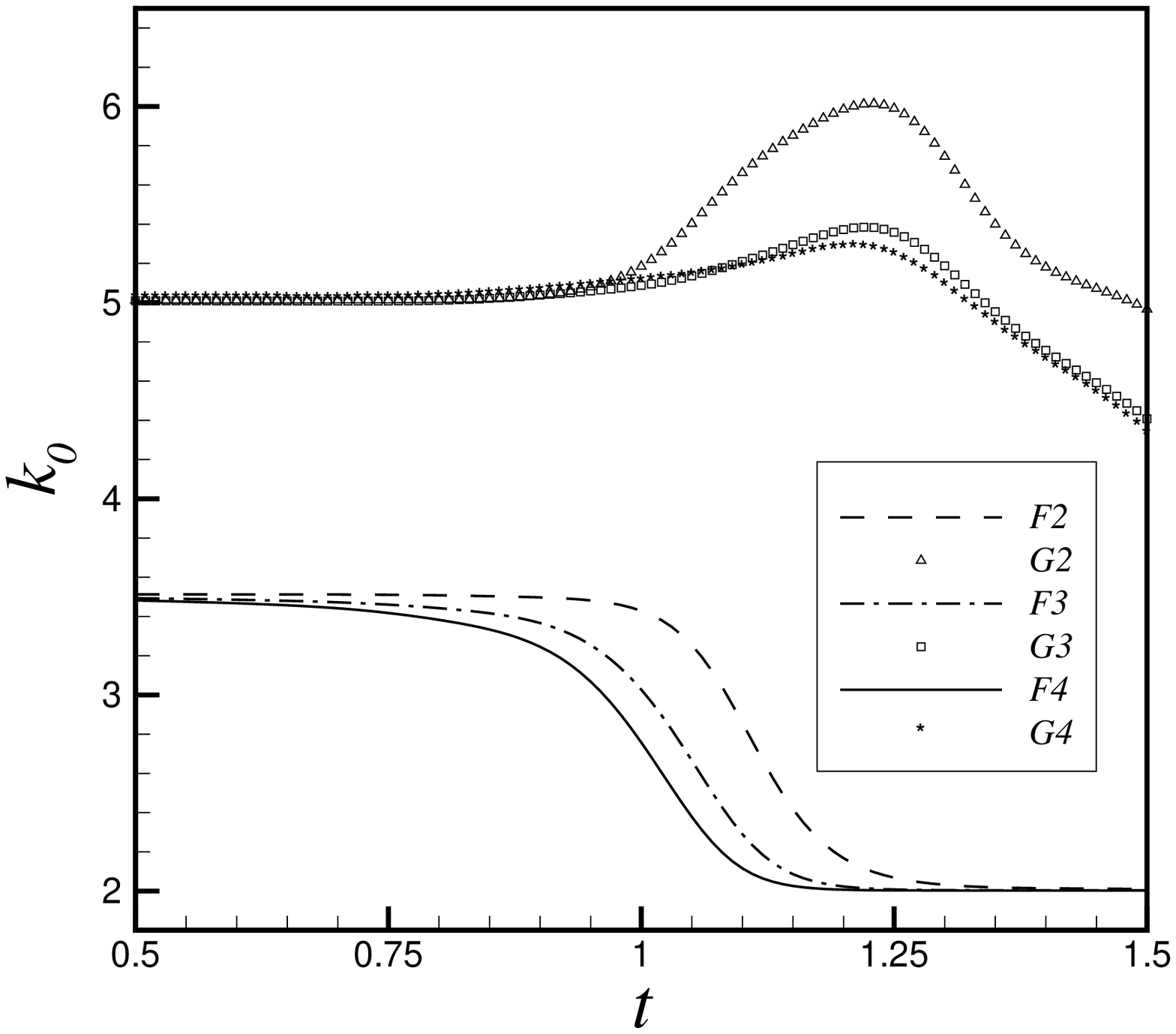, width=6cm}~\epsfig{file=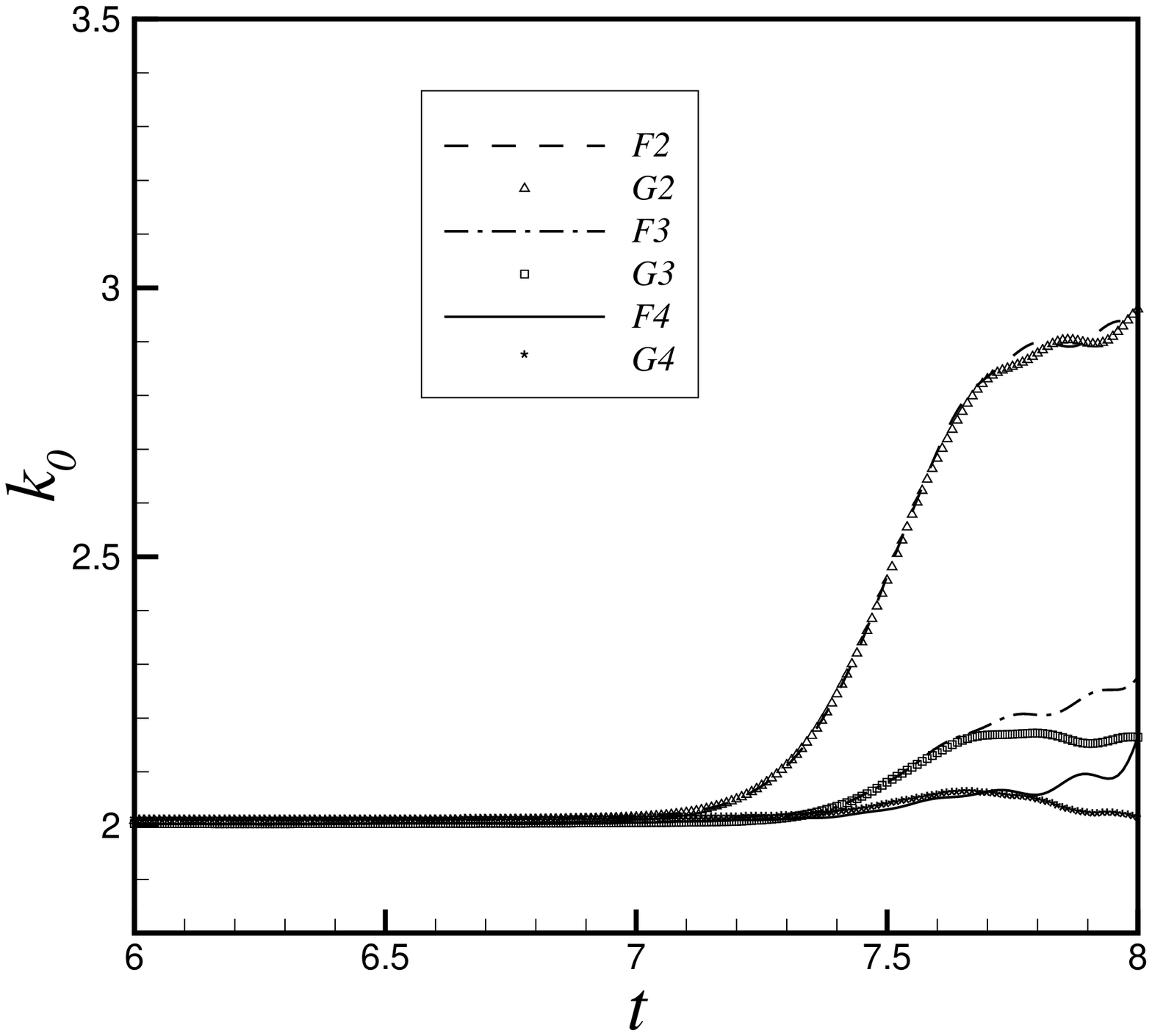, width=6cm}
\caption{Weighted wave numbers $k_0$ as a function of
time. Left: the first phase for $t\in[0.5,~1.5]$; right: the
second phase for $t\in[6,~8]$.}
\end{figure}\end{center}

\vspace{0.5cm}\begin{center} {\bf Table 1} $L_1$-errors $E_1$ and
reflection ratios $r$ for different parameters and grid sizes with
adaptive parameter selection.

{\small
\begin{tabular}{|c|cc|cc|cc|cc|}\hline
 & \multicolumn{2}{|c|}{$E_1$ for $\Delta x=0.1$} &
\multicolumn{2}{|c|}{$E_1$ for $\Delta x=0.05$} &
\multicolumn{2}{|c|}{$r$ for $\Delta x=0.1$} &
\multicolumn{2}{|c|}{$r$ for $\Delta x=0.05$}\\ \hline $p$&
Fourier &
Gabor  & Fourier  & Gabor  & Fourier  & Gabor & Fourier  & Gabor\\
\hline
 1 & 1.62d-2  & 1.79d-2  & 3.28d-2  & 3.47d-2  & 7.43d-3  & 7.63d-3  & 2.87d-2  & 2.94d-2\\ 
 2 & 5.12d-3  & 2.72d-3  & 4.51d-3  & 2.62d-3  & 3.79d-4  & 1.65d-4  & 3.23d-4  & 1.70d-4\\
 3 & 5.29d-3  & 1.93d-3  & 4.69d-3  & 1.54d-3  & 3.66d-4  & 7.32d-5  & 2.99d-4  & 4.48d-5\\
 4 & 5.27d-3  & 1.93d-3  & 5.01d-3  & 1.56d-3  & 3.76d-4  & 7.14d-5  & 3.45d-4  & 4.21d-5\\ 
 5 & 5.07d-3  & 1.95d-3  & 5.07d-3  & 1.59d-3  & 3.70d-4  & 7.23d-5  & 3.70d-4  & 4.49d-5\\ \hline
\end{tabular}}\end{center}

\vspace{0.5cm}
\begin{center}
{\bf Table 2} $L_1$-errors $E_1$ and reflection ratios $r$ for different parameters and grid sizes without adaptivity.

{\small
\begin{tabular}{|c|cc|cc|}\hline
 & \multicolumn{2}{|c|}{$E_1$  } & \multicolumn{2}{|c|}{$r$} \\ \hline
$k_0$& $\Delta x=0.1$ & $0.05$  & $0.1$  & $0.05$ \\ \hline
 2 & 3.22d-3  &  3.26d-3  & 2.00d-4  & 1.73d-4 \\
3.5& 5.33d-3  &  4.98d-3  & 8.58d-4  & 7.89d-4 \\
 5 & 1.26d-2  &  1.23d-2  & 4.81d-3  & 4.60d-3 \\ \hline
\end{tabular}}\end{center}

\vspace{0.5cm}\begin{center}
{\bf Table 3} $L_1$-errors $E_1$ and reflection ratios $r$ for different
window lengths determined adaptively as $b=\beta k_0$.  $p=4.$

{\small
\begin{tabular}{|c|cc|cc|}\hline
 & \multicolumn{2}{|c|}{$E_1$  } & \multicolumn{2}{|c|}{$r$} \\ \hline
$\beta$& $\Delta x=0.1$ & $0.05$  & $0.1$  & $0.05$ \\ \hline
0.5& 7.71d-3  & 2.10e-3   &  1.73e-3 & 1.23e-4  \\
 1 & 1.91d-3  & 1.57d-3   & 7.35d-5  & 4.37d-5 \\
 2 & 1.92d-3  & 1.55d-3   & 6.93d-5  & 4.10d-5 \\
 3 & 1.92d-3  & 1.53d-3   & 7.02d-5  & 3.98d-5 \\
 4 & 1.93d-3  & 1.54d-3   & 7.09d-5  & 4.02d-5 \\ \hline
\end{tabular}}\end{center}
\vspace{0.5cm}

{\bf Example 2.} We then consider a nonlinear wave with repulsive
interaction ($g=2$ in Eq. (\ref{oned})). The initial data and
potential function are taken to be Gaussian pulses
\begin{equation}\psi(x,0)=\mathrm{e}^{-0.1(x-x_0)^2} ~\mathrm{and}~ V(x)=\mathrm{e}^{-0.5(x-x_0)^2}, \end{equation}
with $x_0=15$. This has been an example in \cite{XH:PRE:06} used to
model expansion of a Bose-Einstein condensate which is composed of
waves with different group velocities. The frequency context at the
boundaries is depending on the temporal evolution. In
\cite{XH:PRE:06}, the authors obtained the results under different
wave-number parameters which are independent of time $t$. It was
illustrated that a very bad result appeared if we cannot choose a
suitable $k_0$. Therefore, it is necessary to capture this parameter
adaptively in order to minimize the nonphysical reflection.

In the calculation, $L=30$, $\Delta x=0.1$, and $\Delta t=0.01$ are
chosen. The numerical results with the same mesh sizes by using the
proposed ABC in a large domain $[-15,45]$ are taken to be a
reference solution which is regarded as the ``exact" solution, since
the analytic solution is unknown. Fig. 2 shows the motion of the
wave with the ABCs at time $t=4$ and 6, in which we take $p=4$ and
the window length of Gabor transform $b=L/4$. It is illustrated that
the reflected wave is very small when the waves hit the boundaries
under our adaptive parameter selection strategy. We also show the
wave-number parameters at both boundaries as functions of time in
Fig. 3, in which we see that the wave numbers decay with time after
the waves reach artificial boundaries.

\begin{center}\begin{figure}
 \epsfig{file=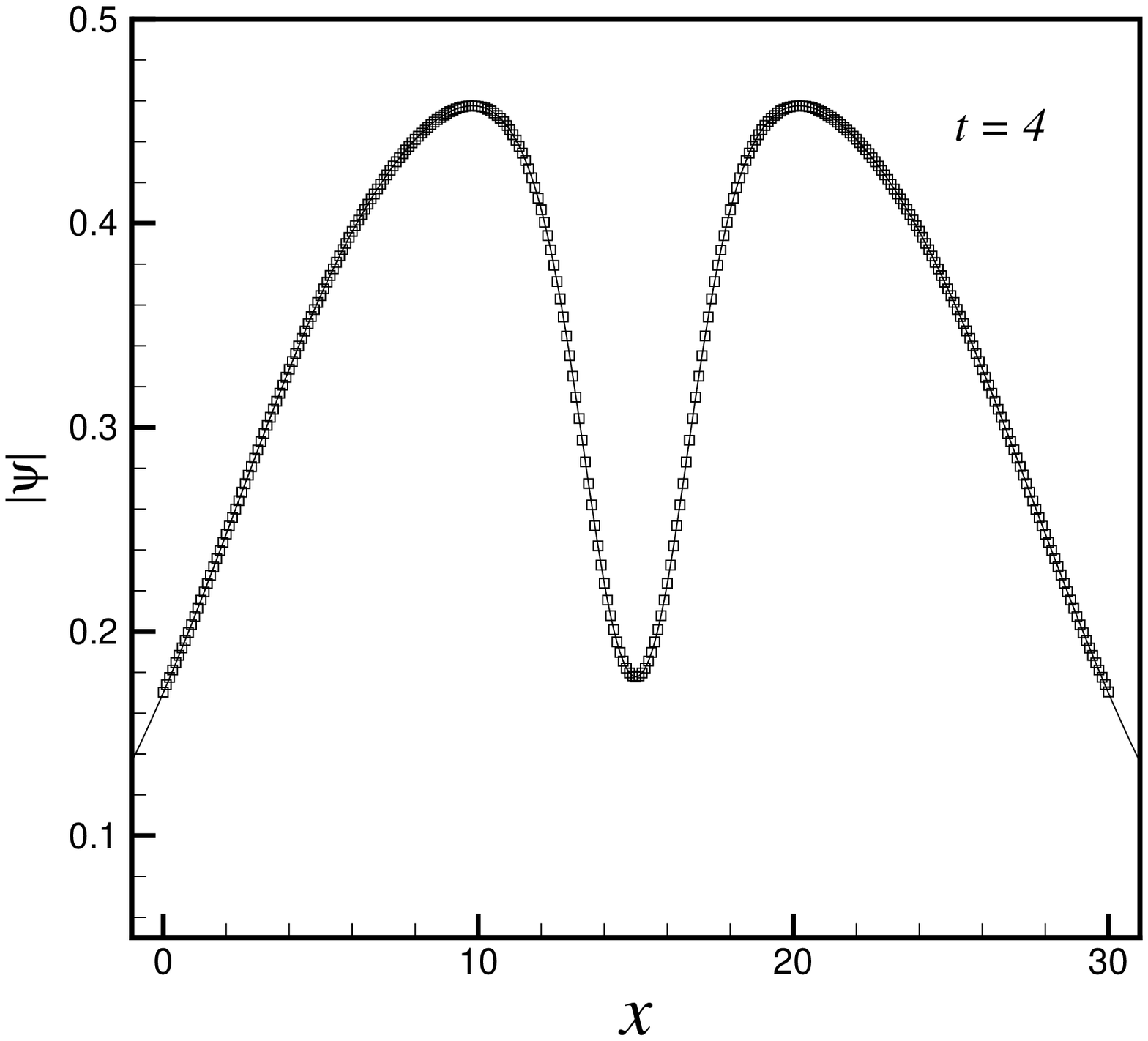, width=6cm}~\epsfig{file=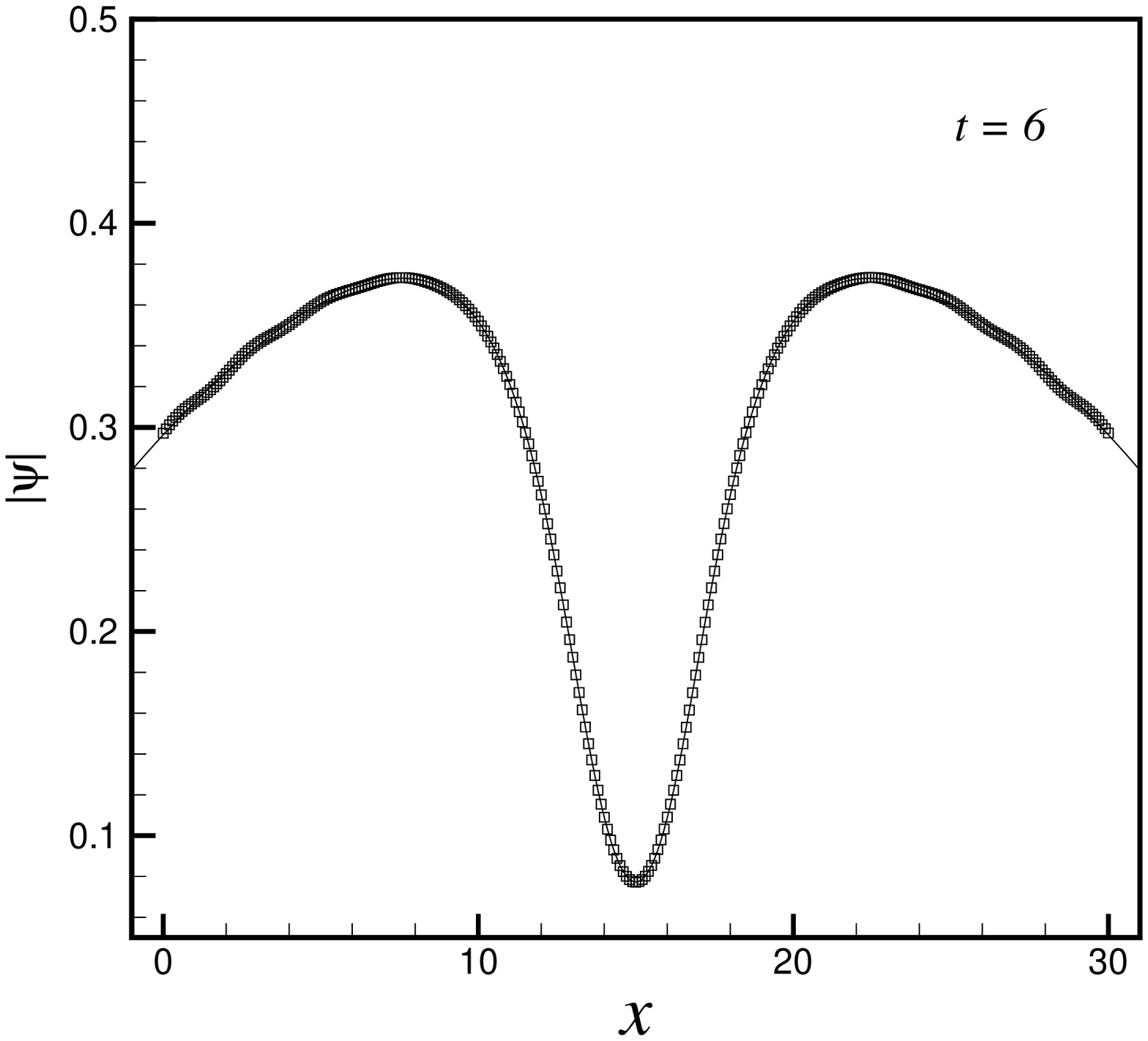, width=6cm}
\caption{The $|\psi|$ solutions at time $t=4$ and 6.}
\end{figure}\end{center}

\begin{center}\begin{figure}
 \epsfig{file=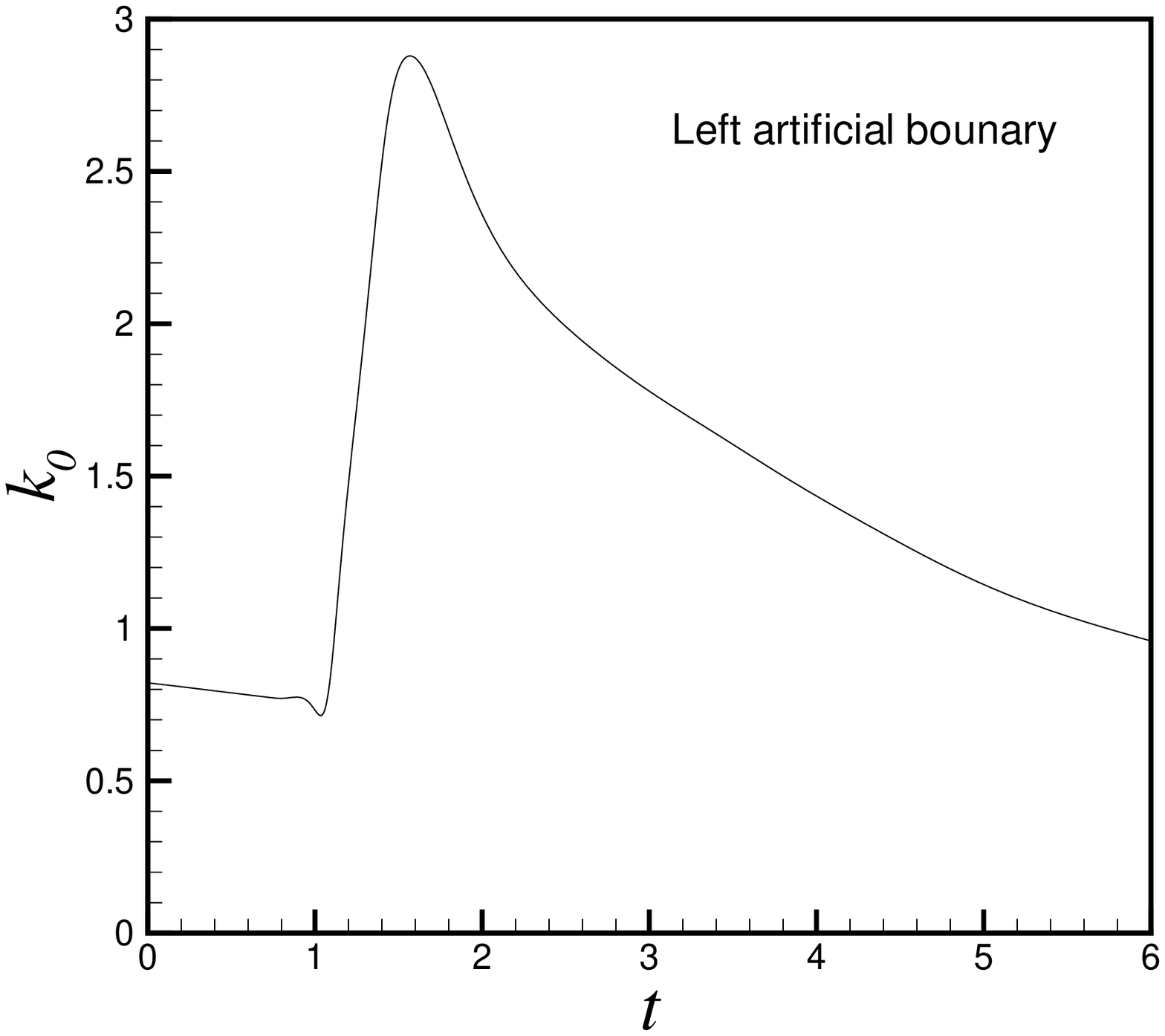, width=6cm}~\epsfig{file=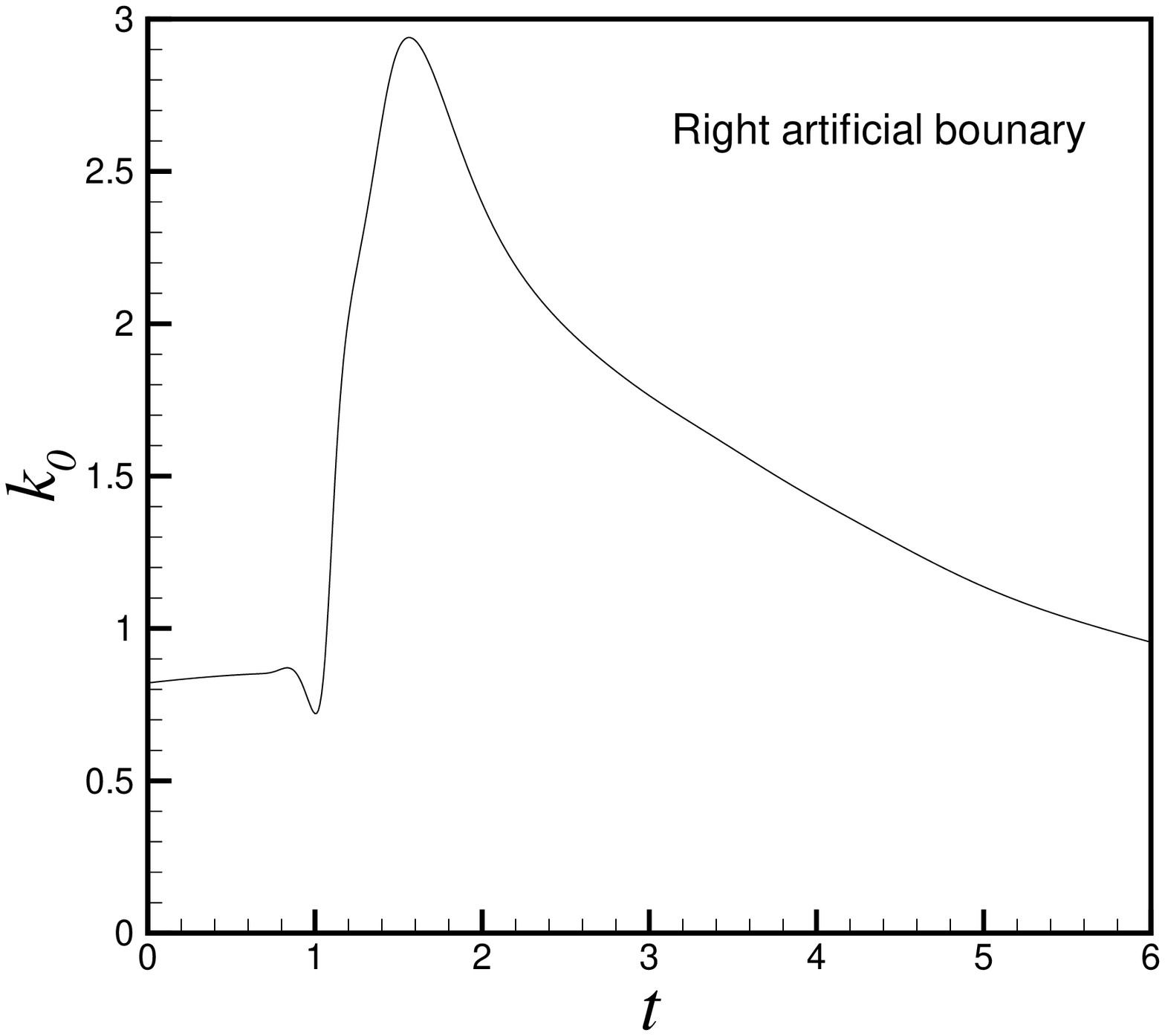, width=6cm}
\caption{ The wave numbers $k_0$ as functions of time $t$.}
\end{figure}\end{center}

{\bf Example 3.} This is a two-dimensional example for Eq.
(\ref{twod}) with cubic nonlinearity $f(|\psi|^2)=-|\psi|^2$ in
homogeneous media; i.e., the potential $V\equiv 0$. We consider the
temporal evolution of an initial packet of the wave centered at
$(x,~y)=(5,~5)$
\begin{equation}
\psi(x,y,0)=\sqrt{2}e^{(x-5)^2+(y-5)^2}e^{2i(x+y-10)}.
\end{equation}
The wave packet moves along the northeast direction and impinges on
the artificial boundaries $\Gamma_e$ and $\Gamma_n$. At the same
time, the amplitude of the wave packet deceases with time due to the
expansion effect. In the calculations, we set the computational
domain be $[0,L]^2$ for $L=10$.  We also set $p=4$ and the window
length of Gabor tranform $b=L/4$. We show numerical solutions of
$|\psi|$ at time $t=0.5,~1,~1.5$ and 2 in Fig. 4 for $h=\Delta
x=\Delta y=0.05$ and $\Delta t=h^2$. We see the wave can be well
absorbed with only minor reflections. In order to see the errors, we
take the numerical result in a large domain $[0,20]^2$ with $h=0.05$
to be a reference solution. In Fig. 5, we plot the temporal
evolution of the solutions and their errors for different mesh sizes at positions
$(x,y)=(10,~10)$ and (10,~5). These results illustrate that the
discussed method can also works well for the two-dimensional
problem.

\begin{center}\begin{figure}
\epsfig{file=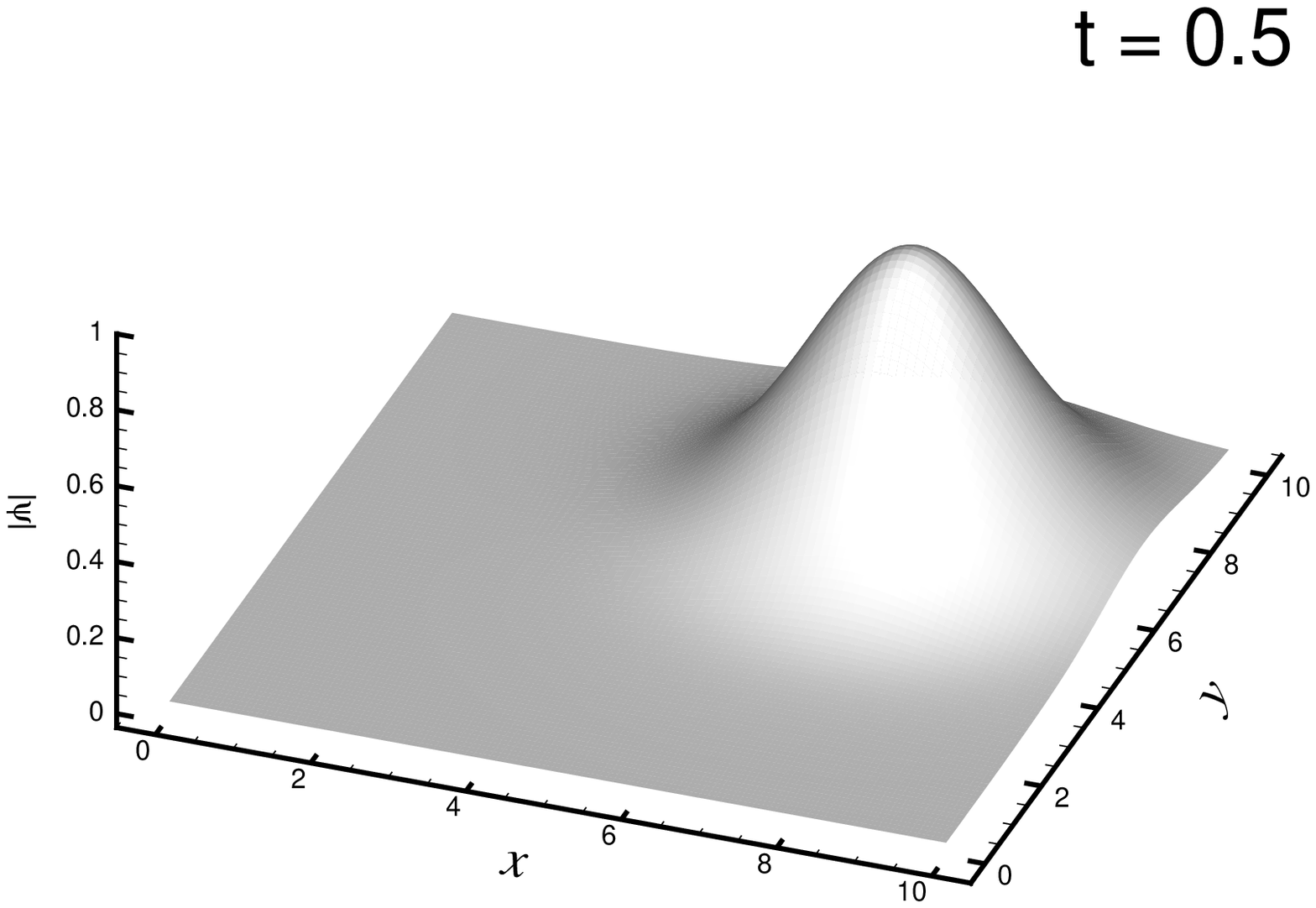, width=7cm}~\epsfig{file=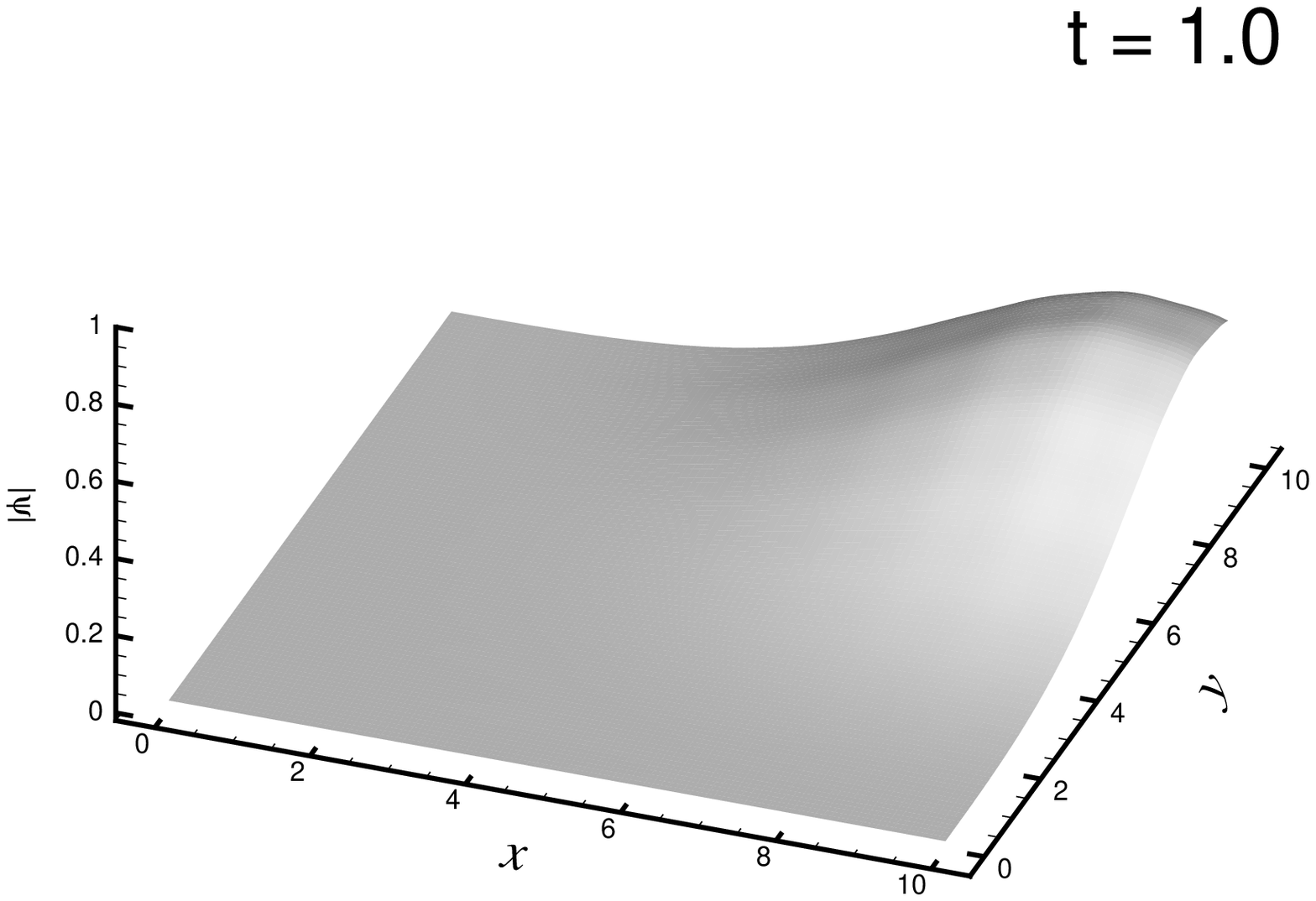, width=7cm}\\
\epsfig{file=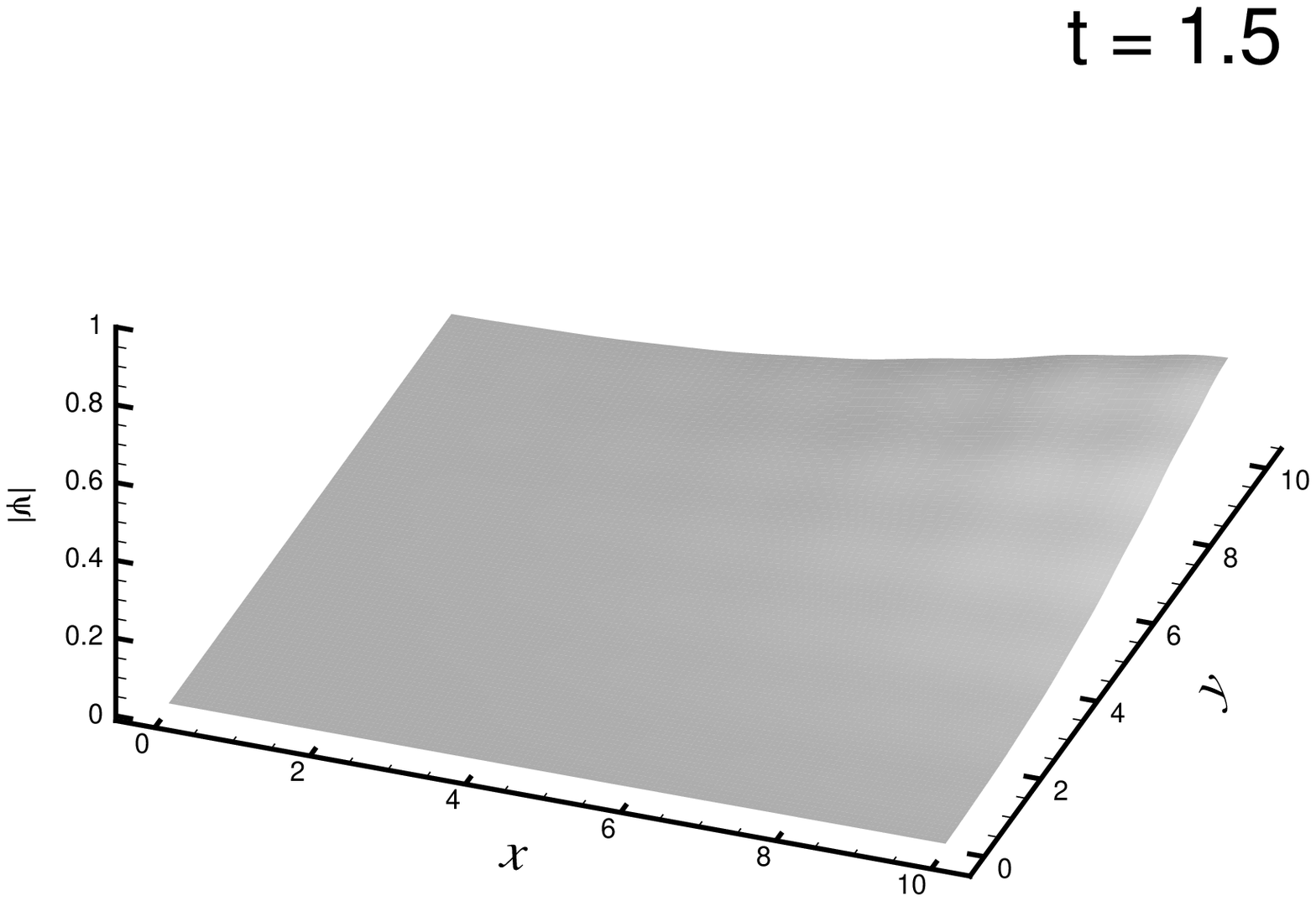, width=7cm}~\epsfig{file=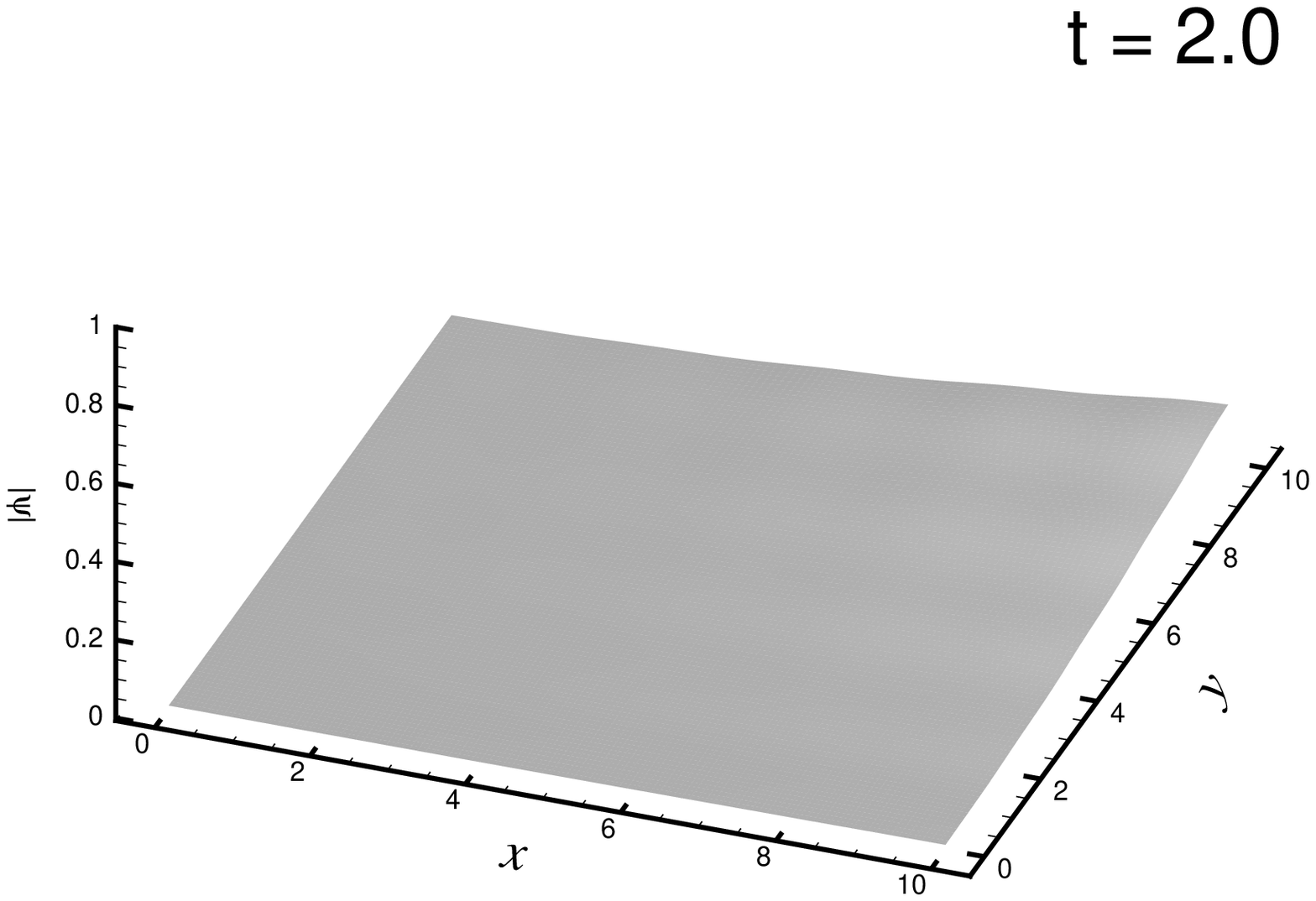, width=7cm}
\caption{Numerical solutions of $|\psi|$ under the mesh
size $\Delta x=\Delta y=0.05$ and $\Delta t=0.0025$ at time
$t=0.5,~1,~1.5$ and 2. }
\end{figure}\end{center}

\begin{center}\begin{figure}
 \epsfig{file=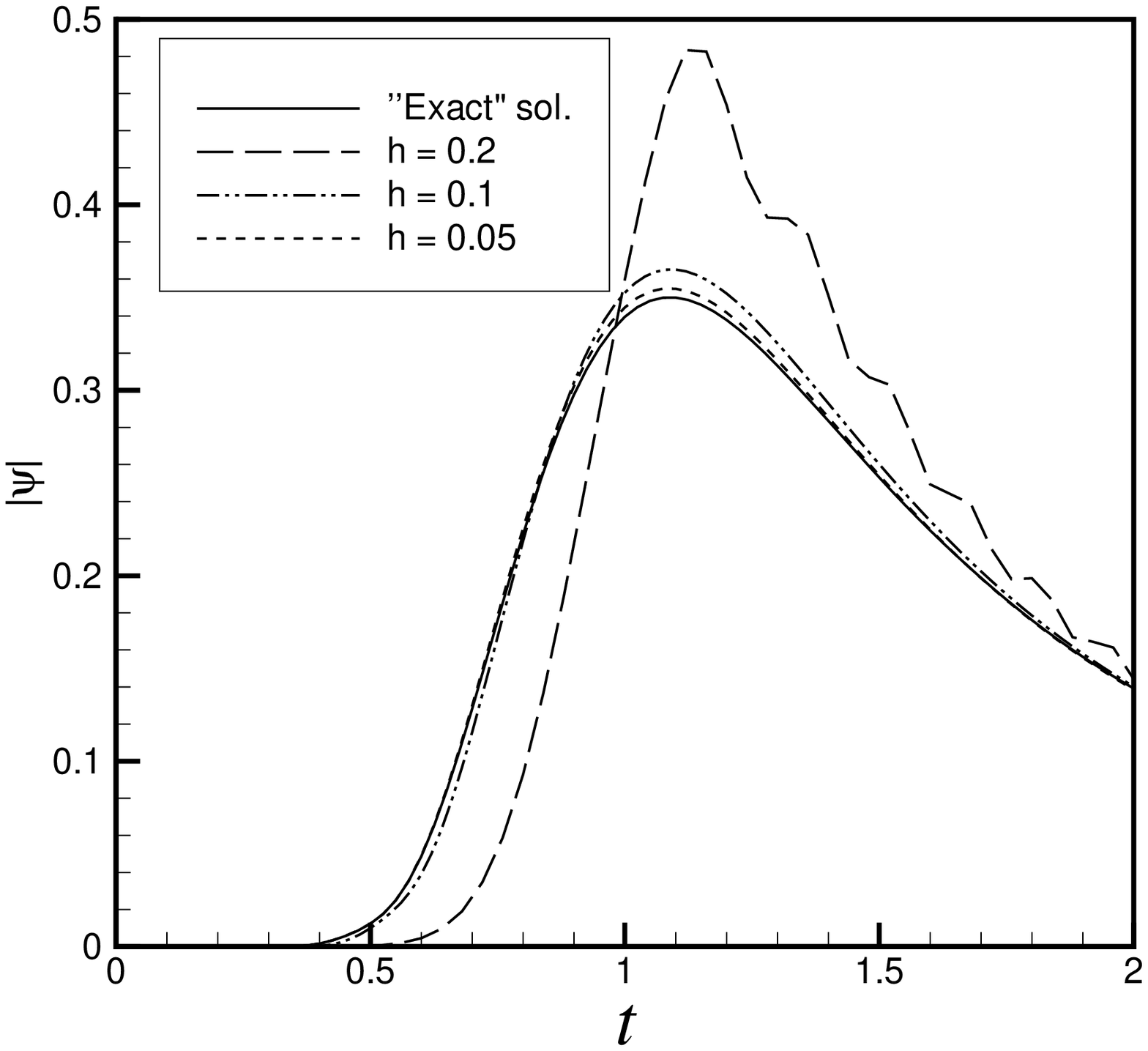, width=6cm}~\epsfig{file=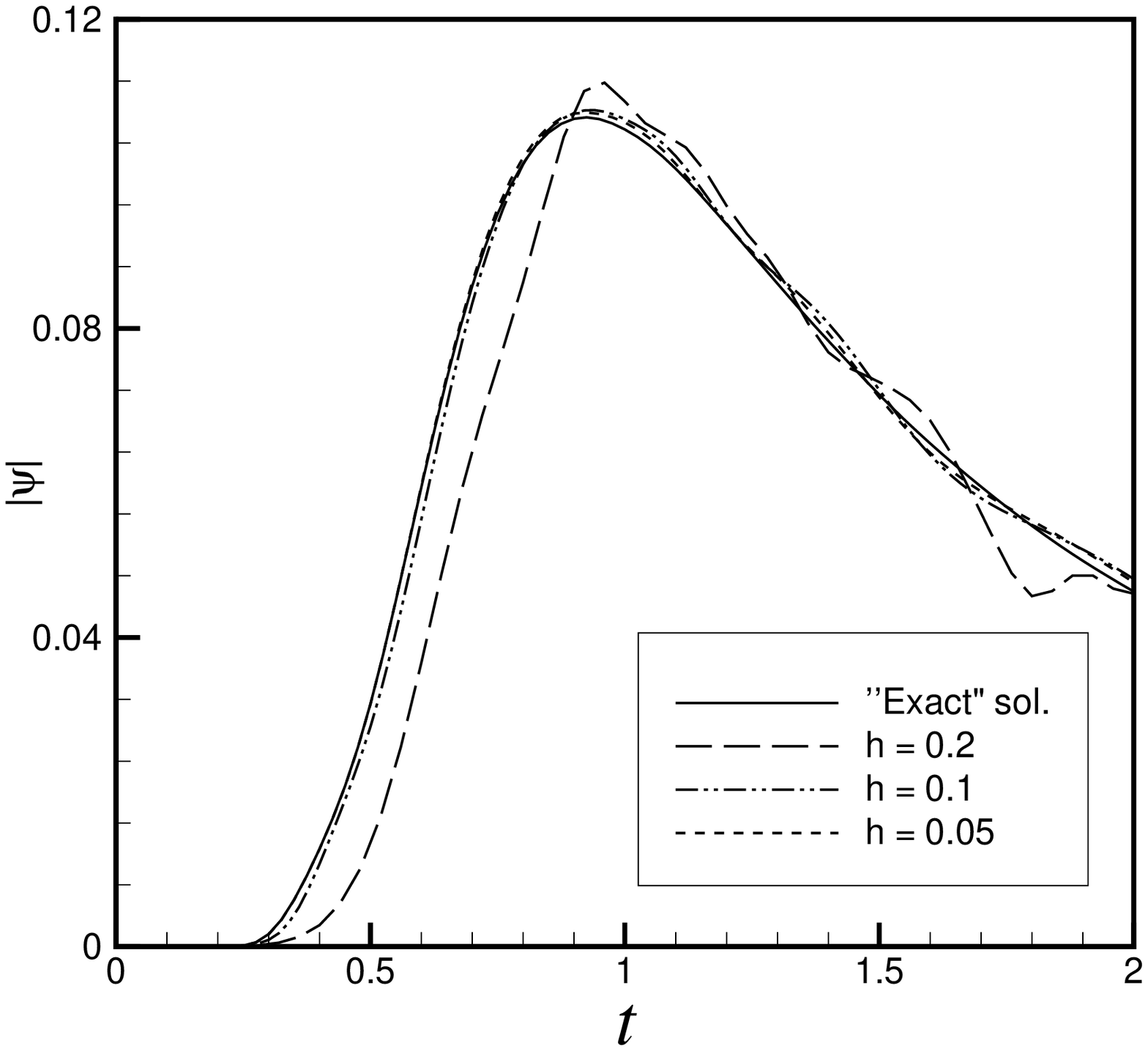, width=6cm}\\
 \epsfig{file=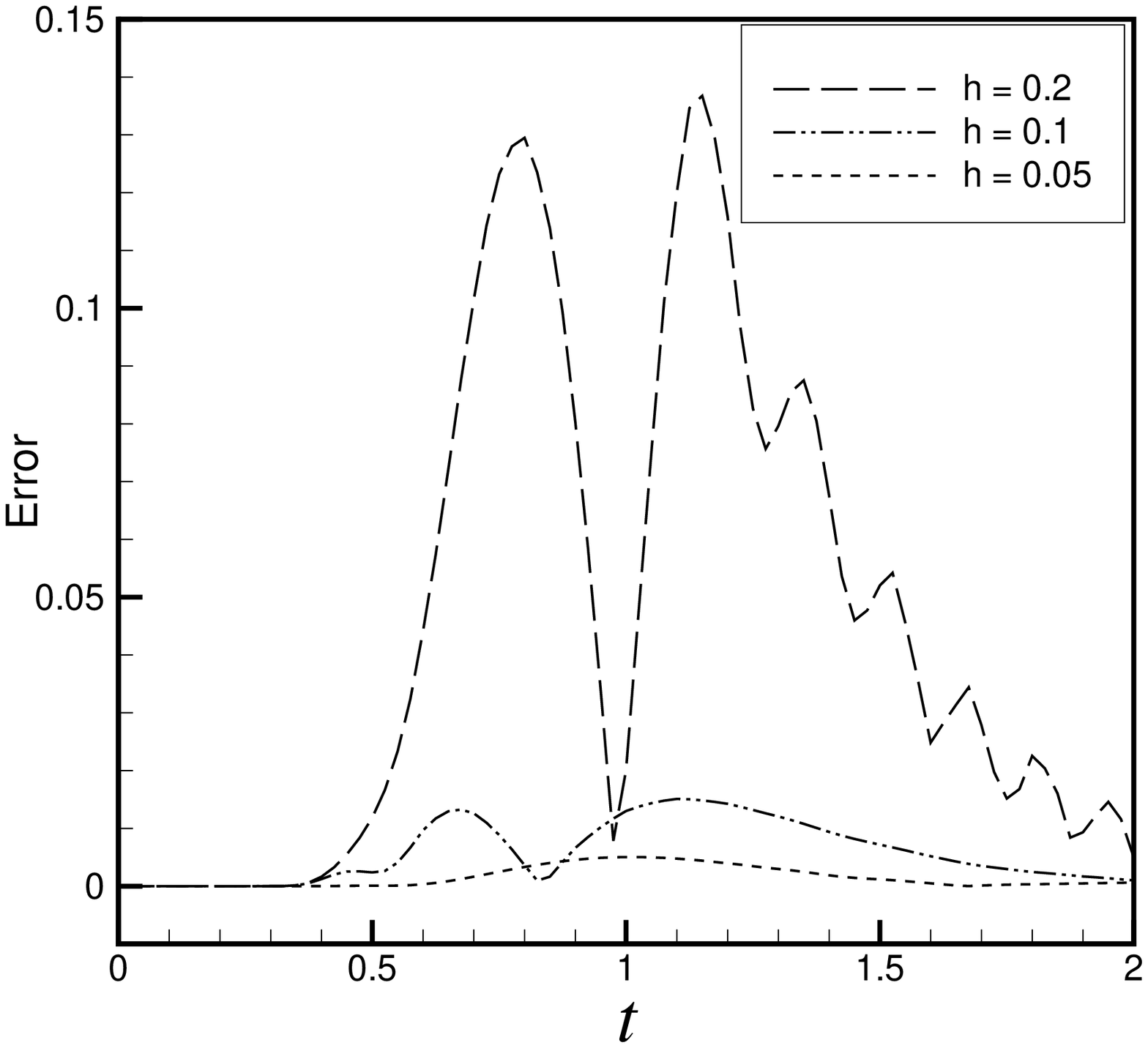, width=6cm}~\epsfig{file=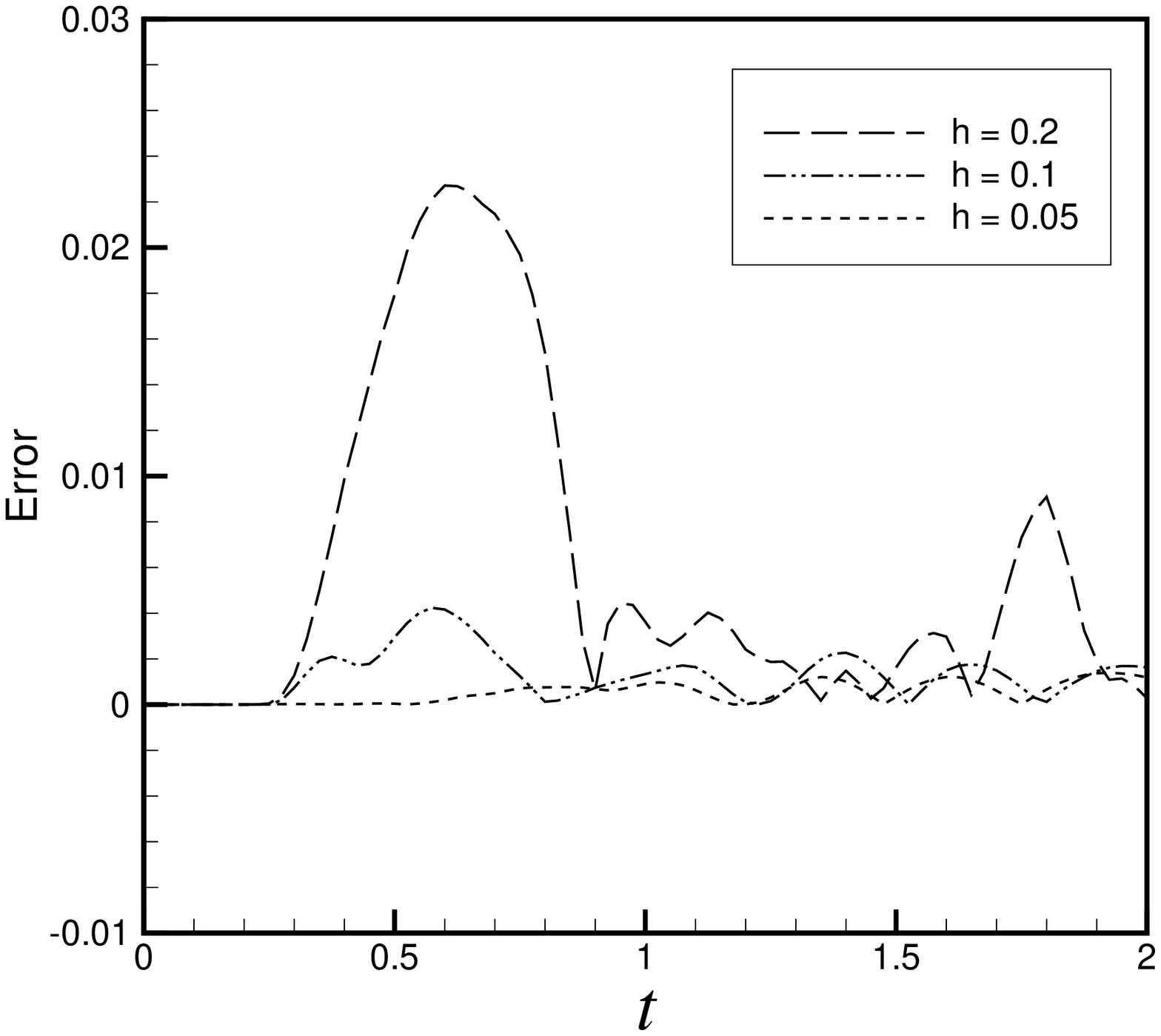, width=6cm}
\caption{Temporal evolutions and errors of $|\psi|$ at positions:
$(x,y)=(10,10)$ (Left), and $(x,y)=(10,5)$ (Right). The ``exact" solution is
computed in a large domain $[0,20]^2$ with the mesh size $h=0.05$}
\end{figure}\end{center}

\section{Concluding remarks}

We develop an efficient adaptive parameter approach for absorbing
boundary conditions of Schr\"{o}dinger-type equations. This approach
is coupled with the local time-splitting method to constitute a
complete procedure for nonlinear problems. We also introduce an
extension to deal with absorbing boundary conditions for
multidimensional nonlinear Schr\"{o}dinger equations. Numerical
examples are performed to show the attractive features of the
approach under consideration. Related further work includes the
stability and error analysis of the proposed approach and further
extension to more complicated initial-boundary value problems.
Another problem is induced by the complexity of nonlinear mechanics.
In some situations, the outgoing waves will return to the interior
domain due to their interactions of nonlinear packets. This open
problem is still unsolved in this paper and we leave for further
consideration.

{\bf Acknowledgments}

This work is supported by the National Natural Science Foundation
of China (Grant No. 10471073), the RGC of Hong Kong and FRG of
Hong Kong Baptist University. The first author thanks Prof. J.
Deng for helpful discussion. The authors also thank the anonymous
referees for many comments and suggestions in the earlier draft of
this paper.


\begin{thebibliography}{10}

\bibitem{SS:Book:99}
C.~Sulem and P.L. Sulem.
\newblock {\em The Nonlinear Schrodinger Equation: Self-Focusing and Wave
  Collapse}.
\newblock Springer, 1999.

\bibitem{Agr:BOOK:01}
G.P. Agrawal.
\newblock {\em Nonlinear fiber optics, 3rd Ed.}
\newblock Academic Press, San Diego, 2001.

\bibitem{GIVOLI:BOOK:92}
D.~Givoli.
\newblock {\em Numerical Methods for Problems in Infinite Domains}.
\newblock Elsevier, Amsterdam, 1992.

\bibitem{Tsy:ANM:98}
S.~V. Tsynkov.
\newblock Numerical solution of problems on unbounded domains. a review.
\newblock {\em Appl. Numer. Math.}, 27:465--532, 1998.

\bibitem{Hag:AN:99}
T.~Hagstrom.
\newblock Radiation boundary conditions for the numerical simulation of waves.
\newblock {\em Acta Numerica}, pages 47--106, 1999.

\bibitem{Han:FPCAM:06}
H.~Han.
\newblock The artificial boundary method -- numerical solutions of partial
  differential equations on unbounded domains.
\newblock In T.~Li and P.~Zhang, editors, {\em Frontiers and Prospects of
  Contemporary Applied Mathematics}, pages 33--58. Higher Education Press and
  World Scientific, 2006.

\bibitem{EM:MC:77}
B.~Engquist and A.~Majda.
\newblock Absorbing boundary conditions for the numerical simulation of waves.
\newblock {\em Math. Comput.}, 31:629--651, 1977.

\bibitem{Higdon:MC:86}
R.L. Higdon.
\newblock Absorbing boundary conditions for difference approximations to the
  multi-dimensional wave equation.
\newblock {\em Math. Comput.}, 47:437--459, 1986.

\bibitem{HW:JCM:85}
H.D. Han and X.N. Wu.
\newblock Approximation of infinite boundary condition and its applications to
  finite element methods.
\newblock {\em J. Comput. Math.}, 3:179--192, 1985.

\bibitem{YU:JCM:85}
D.H. Yu.
\newblock Approximation of boundary conditions at infinity for a harmonic
  equation.
\newblock {\em J. Comput. Math.}, 3:219--227, 1985.

\bibitem{HR:NM:95}
L.~Halpern and J.~Rauch.
\newblock Absorbing boundary conditions for diffusion equations.
\newblock {\em Numer. Math.}, 71:185--224, 1995.

\bibitem{HH:CMA:02-2}
H.~Han and Z.~Huang.
\newblock Exact and approximating boundary conditions for the parabolic
  problems on unbounded domains.
\newblock {\em Comput. Math. Appl.}, 44:655--666, 2002.

\bibitem{SY:JCP:97}
F.~Schmidt and D.~Yevick.
\newblock Discrete transparent boundary conditions for {Schr\"{o}dinger-type
  equations}.
\newblock {\em J. Comput. Phys.}, 134:96--107, 1997.

\bibitem{Arno:VSLID:98}
A.~Arnold.
\newblock Numerically absorbing boundary conditions for quantum evolution
  equations.
\newblock {\em VSLI Design}, 6:313--319, 1998.

\bibitem{AES:CMS:03}
A.~Arnold, M.~Ehrhardt, and I.~Sofronov.
\newblock Discrete transparent boundary conditions for the {Schr\"{o}dinger}
  equation: Fast calculation, approximation, and stability.
\newblock {\em Commun. Math. Sci.}, 1:501--556, 2003.

\bibitem{AB:JCP:03}
X.~Antoine and C.~Besse.
\newblock Unconditionally stable discretization schemes of non-reflecting
  boundary conditions for the one-dimensional {Schr\"{o}dinger} equation.
\newblock {\em J. Comput. Phys.}, 188:157--175, 2003.

\bibitem{HH:CMS:04}
H.~Han and Z.~Huang.
\newblock Exact artificial boundary conditions for {Schr\"{o}dinger equation in
  $R^2$}.
\newblock {\em Commun. Math. Sci.}, 2:79--94, 2004.

\bibitem{SW:JCP:06}
Z.~Z. Sun and X.~Wu.
\newblock The stability and convergence of a difference scheme for the
  {Schr\"{o}dinger} equation on an infinite domain by using artificial boundary
  conditions.
\newblock {\em J. Comput. Phys.}, 214:209--223, 2006.

\bibitem{HJW:CMA:05}
H.~Han, J.~Jin, and X.~Wu.
\newblock A finite-difference method for the one-dimensional time-dependent
  {Schr\"{o}dinger} equation on unbounded domain.
\newblock {\em Comput. Math. Appl.}, 50:1345--1362, 2005.

\bibitem{JG:CMA:04}
S.~D. Jiang and L.~Greengard.
\newblock Fast evalution of nonreflecting boundary conditions for the
  {Schr\"{o}dinger} equation in one dimension.
\newblock {\em Comput. Math. Appl.}, 47:955--966, 2004.

\bibitem{Shi:PRB:91}
T.~Shibata.
\newblock Absorbing boundary conditions for the finite-difference time-domain
  calculation of the one dimensional {Schr\"{o}dinger equation}.
\newblock {\em Phys. Rev. B}, 43:6760, 1991.

\bibitem{Kus:PRB:92}
J.-P. Kuska.
\newblock Absorbing boundary conditions for the {Schr\"{o}dinger equation} on
  finite intervals.
\newblock {\em Phys. Rev. B}, 46:5000, 1992.

\bibitem{Di:NFAO:97}
L.~{Di Menza}.
\newblock Transparent and absorbing boundary conditions for the
  {Schr\"{o}dinger} equation in a bounded domain.
\newblock {\em Numer. Funct. Anal. Optim.}, 18:759--775, 1997.

\bibitem{FJ:SISC:99}
T.~Fevens and H.~Jiang.
\newblock Absorbing boundary conditions for the {Schr\"{o}dinger equation}.
\newblock {\em SIAM J. Sci. Comput.}, 21:255--282, 1999.

\bibitem{AR:SINA:02}
I.~Alonso-Mallo and N.~Reguera.
\newblock Weak ill-posedness of spatial discretizations of absorbing boundary
  conditions for {Schr\"{o}dinger-type equations}.
\newblock {\em SIAM J. Numer. Anal.}, 40:134--158, 2002.

\bibitem{Sze:SINA:04}
J.~Szeftel.
\newblock Design of absorbing boundary conditions for {Schr\"{o}dinger}
  equations in {$R^d$}.
\newblock {\em SIAM J. Numer. Anal.}, 42:1527--1551, 2004.

\bibitem{HK:MC:87}
T.~Hagstrom and H.B. Keller.
\newblock Asymptotic boundary conditions and numerical methods for nonlinear
  elliptic problems on unbounded domains.
\newblock {\em Math. Comput.}, 48:449--470, 1987.

\bibitem{HWX:JCM:06}
H.D. Han, X.N. Wu, and Z.L. Xu.
\newblock Artificial boundary method for {Burgers'} equation using nonlinear
  boundary conditions.
\newblock {\em J. Comput. Math.}, 24:295--304, 2006.

\bibitem{XHW:CCP:06}
Z.~Xu, H.~Han, and X.~Wu.
\newblock Numerical method for the deterministic {Kardar-Parisi-Zhang} equation
  in unbounded domains.
\newblock {\em Commun. Comput. Phys.}, 1:481--495, 2006.

\bibitem{Zhe:JCP:06}
C.~Zheng.
\newblock Exact nonreflecting boundary conditions for one-dimensional cubic
  nonlinear {Schr\"{o}dinger} equations.
\newblock {\em J. Comput. Phys.}, 215:552--565, 2006.

\bibitem{ABD:SINA:06}
X.~Antoine, C.~Besse, and S.~Descombes.
\newblock Artificial boundary conditions for one-dimensional cubic nonlinear
  {Schr\"{o}dinger} equations.
\newblock {\em SIAM J. Numer. Anal.}, 43:2272--2293, 2006.

\bibitem{Sze:CMAME:06}
J.~Szeftel.
\newblock Absorbing boundary conditions for nonlinear scalar partial
  differential equations.
\newblock {\em Comput. Methods Appl. Mech. Engrg.}, 195:3760--3775, 2006.

\bibitem{Sze:NM:06}
J.~Szeftel.
\newblock Absorbing boundary conditions for one-dimensional nonlinear
  {Schr\"{o}dinger} equations.
\newblock {\em Numer. Math.}, 104:103--127, 2006.

\bibitem{FL:JOB:05}
C.~Farrell and U.~Leonhardt.
\newblock The perfectly matched layer in numerical simulations of nonlinear and
  matter waves.
\newblock {\em J. Opt. B: Quantum Semiclass. Opt.}, 7:1--4, 2005.

\bibitem{XH:PRE:06}
Z.~Xu and H.~Han.
\newblock Absorbing boundary conditions for nonlinear {Schr\"{o}dinger
  equations}.
\newblock {\em Phys. Rev. E}, 74:037704, 2006.

\bibitem{Gabor:46}
D.~Gabor.
\newblock Theory of communication.
\newblock {\em J. Inst. Elect. Eng. (London)}, 93:429--457, 1946.

\bibitem{SSP:CPAM:84}
P.L. Sulem, C.~Sulem, and A.~Patera.
\newblock Numerical simulation of singular solutions to the two-dimensional
  cubic {Schr\"{o}dinger} equation.
\newblock {\em Commum. Pure Appl. Math.}, 37:755--778, 1984.

\bibitem{CJS:JCP:99}
Q.~Chang, E.~Jia, and W.~Sun.
\newblock Difference schemes for solving the generalized nonlinear
  {Schr\"{o}dinger} equation.
\newblock {\em J. Comput. Phys.}, 148:397--415, 1999.

\bibitem{BJM:SISC:03}
W.~Bao, S.~Jin, and P.A. Markowich.
\newblock Numerical study of time-splitting spectral discretizations of
  nonlinear {Schr\"{o}dinger} equations in the semiclassical regimes.
\newblock {\em SIAM J. Sci. Comput.}, 25:27--64, 2003.

\bibitem{Str:SINA:68}
G.~Strang.
\newblock On the constrction and comparison of difference schemes.
\newblock {\em SIAM J. Numer. Anal.}, 5:506--517, 1968.

\end{thebibliography}

\end{document}